\documentclass[12pt,a4paper]{amsart}
\usepackage[height=715pt,width=420pt,centering]{geometry}
\usepackage{mathtools}          
\usepackage{amssymb}
\usepackage{xspace}
\usepackage{tikz}

\usepackage[plainpages=false,colorlinks=false]{hyperref}

\theoremstyle{plain}
\newtheorem{theorem}{Theorem}[section]
\newtheorem{proposition}[theorem]{Proposition}
\newtheorem{lemma}[theorem]{Lemma}
\newtheorem{corollary}[theorem]{Corollary}
\theoremstyle{definition}
\newtheorem{definition}[theorem]{Definition}

\theoremstyle{remark}
\newtheorem{remark}[theorem]{Remark}

\tikzstyle{vertex}=[circle, draw, fill=black!0, inner sep=1pt, minimum width=3pt]
\allowdisplaybreaks[1]          

\DeclareMathOperator{\Hom}{Hom}
\DeclareMathOperator{\Irr}{Irr}

\let\phi=\varphi

\let\join=\vee
\let\cross=\times
\let\normal=\lhd

\let\Direct=\bigoplus
\let\intersect=\cap
\let\Intersect=\bigcap
\let\tensor=\otimes
\let\isomorph=\cong

\newcommand\idx[2]{\abs{#1:#2}}

\newcommand\class{\mathcal{K}}
\newcommand\classsum{\hat{\mathcal{K}}}
\newcommand\classS[1]{\abs{\class_{#1}}}
\newcommand\phan{{\vphantom{1}}}
\newcommand\ie{\emph{i.e.,}\@\xspace} 
\newcommand\eg{\emph{e.g.,}\@\xspace} 
\newcommand\inv[1]{\frac{1}{#1}}
\newcommand\abs[1]{\left\lvert #1 \right\rvert}
\newcommand\angles[1]{ \left\langle #1 \right\rangle }
\newcommand\parens[1]{\left(#1\right)}

\newcommand\defn[1]{#1\xspace}
\newcommand\presentation[2]{\,\left\langle #1 \;\middle\vert\;#2
    \vphantom{{#1#2}'} \right\rangle\,}

\title{Transposable Character Tables, Dual Groups}
\author{Ivan Andrus}
\address{%
  Department of Mathematics and its Applications\\
  Central European University\\
  Nador utca 9\\
  1051 Budapest, Hungary%
}
\email{Andrus\textunderscore Ivan@ceu-budapest.edu}

\author{P\'al Heged\H us}
\thanks{The research of the second author was partly supported by OTKA (K 84233).}
\address{%
  Department of Mathematics and its Applications\\
  Central European University\\
  Nador utca 9\\
  1051 Budapest, Hungary%
}
\email{HegedusP@ceu.hu}

\author{Tetsuro Okuyama}
\address{%
  Laboratory of Mathematics\\
  Hokkaido University of Education\\
  Asahikawa, Hokkaido 070-0825, Japan%
}
\email{okuyama.tetsuro@a.hokkyodai.ac.jp}


\begin{document}

\begin{abstract}
  One way of expressing the self-duality $A\cong \Hom(A,\mathbb{C})$ of Abelian groups is that their character tables are self-transpose (in a suitable ordering).
  In this paper we extend the duality to some noncommutative groups considering when the character table of a finite group is close to being the transpose of the character table for some other group.
  We find that groups dual to each other have dual normal subgroup lattices. We show that our concept of duality cannot work for non-nilpotent groups and we describe $p$-group examples.
\end{abstract}

\maketitle

\section{Introduction}


Conjugacy classes and irreducible characters of finite groups can be considered ``dual'' in many ways.
Multiplicative structures defined on them have similar properties connected to the fact that in the character table they satisfy similar orthogonality relations.
The number theoretical properties of the sizes of the classes and the degrees of the irreducible characters imply similar consequences.
See for example the survey articles~\cite{Lewis:2008a} and~\cite{Camina:2011} which have a few such comparisons.
For a direct duality see Propositions~\ref{thm:commutator} and~\ref{thm:char-central-series} in this paper.
There are differences, however, as the following illustrate. (In the graphs referred to here $a,b>1$ are connected if $(a,b)>1$.)
\begin{theorem}[{\cite{Bianchi:2007}}]
  If the character degree graph is complete then $G$ is solvable.
\end{theorem}
\begin{theorem}[{\cite{Fisman:1987}}]
  If the group is nonabelian simple then the conjugacy class size graph is complete.
\end{theorem}

The character table provides a duality between the columns and the rows, but in general this does not extend to a duality of one group to another group.
In this paper we explore when there is a direct way of defining the dual group of a group.

One obvious case is when the transpose of a character table $M$ is itself a character table.
The trivial conjugacy class is the only column of the character table consisting of all positive integers, and hence must correspond in the transpose to the trivial character which is always a row of ones.
This then implies that all character degrees are one which is only true for Abelian groups.
In fact, Abelian groups are $\mathbb{Z}$-modules and for $\mathbb{Z}$-modules there is a concept of duality $M\mapsto M^*=\Hom(M,\mathbb{C})$.
That $M\cong M^*$ follows from the structure of finite Abelian groups.
A main motivation for this paper is to generalize this.

One way to extend the notion of duality to non-Abelian groups is the following.

\begin{definition}
  We say that a group $G$ with character table $X$ is \defn{transposable} if there exist non-negative diagonal integer matrices $D$ and $N$ such that $\tilde{X}=(D^{-1}XN)^{T}$ is the character table of some group $G^{T}$.
  We also say that $G^{T}$ is a \defn{transpose group} of $G$.
\end{definition}

\begin{remark}
  In our definition we stipulated that the matrices $D$ and $N$ must be matrices of non-negative integers.
  It is clear that it is sufficient to consider such matrices.
  In fact, $D$ is the diagonal matrix consisting of the character degrees of $G$, and $N$ is of the degrees of $G^T$.
\end{remark}

With this definition we find that, although $G^{T}$ is not unique, its character table is.
We also show that duality of $G$ and $G^{T}$ imply that their normal subgroup lattices are dual to each other, this is another similarity to the behaviour of dual modules. Within this duality the members of the lower central series of $G$ correspond to the members of the upper central series of $G^{T}$ and vice versa.

We prove that our definition can only work for nilpotent groups and is inherited by direct products and direct factors.
So only $p$-groups need be considered.

In the final section we give some examples of non-Abelian transposable $p$-groups.

Note that this concept was first outlined by Eiitchi Bannai~\cite{Bannai:1993} motivated by group association schemes and fusion algebras. His notation and approach is much different from ours.
He defined self-duality for groups and investigated groups of order $64$.
Subsequently, examples of self-dual groups were found by Akihide Hanaki and the third author~\cite{Hanaki:1997}, and some more examples by Akihide Hanaki~\cite{Hanaki:1997a}.
The first was motivated by the Suzuki $2$-groups of Graham Higman, whose character tables have also been studied by Jeffrey M. Riedl~\cite{Riedl:1999} and I. A. Sagirov~\cite{Sagirov:2003}.


\section{Basic Results}


Of immediate concern are the possible values for the new character degrees.
Perhaps unsurprisingly, they are unique and given by the square roots of the conjugacy class sizes. This result is inherent in {\cite{Bannai:1993}}.

\begin{proposition}
  \label{thm:uniq-char-degr}
  Let $X$ be the character table of a finite group $G$.
  Let $D$ be the diagonal matrix with the character degrees of $G$ along the diagonal (in the same order as they appear in $X$).
  Let $N$ be an arbitrary diagonal matrix such that $\tilde{X}=(D^{-1}XN)^{T}=NX^{T}D^{-1}$ is the character table of some group $G^{T}$.
  Then $N^{2}$ is the diagonal matrix with entries equal to the conjugacy class sizes of $G$ (in the same order as they appear in $X$).
\begin{proof}
  If $\tilde{X}$ is the character table of some group then it must satisfy the orthogonality relations, in particular column orthogonality.
  The columns of $\tilde{X}$ are indexed by irreducible characters of $G$, and the rows by its conjugacy classes.
  We indicate by $d_{i}$ the new character degrees \ie the diagonal entries of $N$, and by $g_{i}$ a representative of the $i$th conjugacy class of $G$ (as indexed by $X$).
  Assume that $g_{1}$ is the identity element, hence $d_{1}=1$, and denote by $1_{G}$ the trivial character of $G$.
  Let $n$ be the number of irreducible characters of $G$.

  The new character table $\tilde{X}$ has $d_i\frac{\chi(g_i)}{\chi(1)}$ in row $i$, column $\chi$.
  Column orthogonality of $\tilde{X}$ implies that for $\chi\not=1_{G}$
  \begin{align*}
    0 & = \sum_{i=1}^{n} d_{i} \frac{\chi(g_{i})}{\chi(1)} \overline{d_{i}} = \inv{\chi(1)}\sum_{i=1}^{n} d_{i}^{2} \chi(g_{i}) = \sum_{i=1}^{n} d_{i}^{2} \chi(g_{i}),
  \end{align*}
  which constrains the possible values for the $d_{i}$.
  Letting $k_{i}^{\phan}=d_{i}^{2}$, we obtain a system of $n-1$ linear equations in $n-1$ unknowns.
  From the row orthogonality relations for $X$ we know that $k_{i}^{\phan}=\abs{g_{i}^{G}}$ is a solution.
  Since the equations are simply $n-1$ of the rows of $X$ they are linearly independent, and so the solution is unique.
\end{proof}
\end{proposition}

\begin{corollary}
  If a finite group $G$ is transposable, $\class$ a conjugacy class of $G$, $x\in \class$ and $\chi\in\Irr(G)$  then $\frac{\sqrt{\abs{\class}}\chi(x)}{\chi(1)}$ is an algebraic integer, in particular the conjugacy classes of $G$ must all have square size.
\begin{proof}The expression $\frac{\sqrt{\abs{\class}}\chi(x)}{\chi(1)}$ is a character table entry for the transpose group $G^T$, so it must be an algebraic integer.
\end{proof}
\end{corollary}

Very little seems to be known about groups with only square conjugacy class sizes.

It is often convenient to define the \defn{normalized character table} of $G$ to be the character table with the rows divided by the character degrees. So $G$ and $G^T$ are transposes of each other if and only if their normalized character tables are transposes of each other.

It is easy to see that transposability is closed under direct products, and in fact it is determined entirely by the direct factors.

\begin{proposition}
  \label{thm:factors}
  Suppose that $M$ is the character table of a group $G$ and factors as the Kronecker product of two matrices $X$ and $Y$ such that each first row of $X$ and $Y$ consists of ones and each first column of positive integers.
  Then $G$ is as a non-trivial direct product of groups having character tables $X$ and $Y$.
\begin{proof}
  By assumption of the structure of $M$, we have character degrees $a_{i}=(X)_{i,1}$, $b_{i}=(Y)_{i,1}$ such that $a_{1}=b_{1}=1$,
  \begin{equation*}
    X=
    \begin{pmatrix}
      1      & 1 & \cdots \\
      b_{2}  & * &        \\
      \vdots &   &        \\
    \end{pmatrix},
    \qquad
    Y=
    \begin{pmatrix}
      1      & 1 & \cdots \\
      a_{2}  & * &        \\
      \vdots &   &        \\
    \end{pmatrix},
  \end{equation*}
  and $M$ has the form
  \begin{center}
    \begin{tabular}{c|ccc|ccc|ccc|c}
      \multicolumn{1}{c}{} & \multicolumn{3}{|c|}{$\mathcal{A}$}                                                                                                                                                                                                                     \\
      \cline{1-11}
      $A$                  & $1$                                 & $1$               & $\dots$              & $1$                                 & $1$               & $\dots$              & $1$                                 & $1$               & $\dots$ &                   \\
      \cline{1-1}
                           & $a_{2}$                             & \phantom{$b_{2}$} &                      & $a_{2}$                             & \phantom{$b_{2}$} &                      & $a_{2}$                             & \phantom{$b_{1}$} &         & $B$               \\
                           & $\vdots$                            &                   &                      & $\vdots$                            &                   &                      & $\vdots$                            &                   &                             \\
      \cline{1-11}
      $A$                  & $b_{2}$                             & $b_{2}$           & $\dots$              &                                     &                   &                      &                                     &                   &         &                   \\
      \cline{1-1}
                           & $a_2b_2$                            &                   &                      &                                     &                   &                      &                                     &                   &                             \\
                           &                                     &                   &                      &                                     &                   &                      &                                     &                   &                             \\
      \cline{1-10}
      $A$                  & $b_{3}$                             & $b_{3}$           & $\dots$              &                                     &                   &                      &                                     &                   &         &                   \\
      \cline{1-1}
                           &                                     &                   &                      &                                     &                   &                      &                                     & $\ddots$          &                             \\
                           &                                     &                   &                      &                                     &                   &                      &                                     &                   &         & $\vdots$          \\
      \cline{2-10}
      \multicolumn{1}{c}{} & \multicolumn{1}{|c|}{$\mathcal{B}$} &                   & \multicolumn{1}{c}{} & \multicolumn{1}{|c|}{$\mathcal{B}$} &                   & \multicolumn{1}{c}{} & \multicolumn{1}{|c|}{$\mathcal{B}$} &                   & \multicolumn{1}{c}{$\dots$} \\
    \end{tabular}
  \end{center}
  Note that the upper left block is the matrix $Y$, and the upper left entries of each block form the matrix $X$.

  Let $\mathcal{A}$ be the subset of conjugacy classes consisting of the first column of blocks.
  Also, let $A$ be the subset of rows consisting of the first row of each block as marked.
  We claim that the intersection of the kernels of these characters
  \[N=\Intersect_{\chi\in A}\ker\chi=\ker\Direct_{\chi\in A}\chi\] %
  is the normal subgroup containing exactly the conjugacy classes in $\mathcal{A}$.

  It is clear that $N\supset\mathcal{A}$, since when restricted to $\mathcal{A}$ the characters are all integer multiples of the trivial character.
  Suppose $N\supsetneq\mathcal{A}$, then some column of $M$ outside of $\mathcal{A}$ is such that the first entry in the $i$th block is given by $b_{i}$ for all $i$.
  Since these entries are $1\cdot b_{i,j}$ it must be that $b_{i,j}=b_{i}$ for $j>1$ and all $i$.
  This is impossible since the columns of $M$ are linearly independent and hence the columns of $X$ must be as well.
  Thus $N=\mathcal{A}$.

  Likewise, let $B$ be the first block of rows and $\mathcal{B}$ the first column in each block.
  The intersection of kernels $H=\Intersect_{\chi\in B}\ker\chi=\ker\Direct_{\chi\in B}\chi$ corresponds to the subset $\mathcal{B}$ of conjugacy classes.
  In particular $N$ and $H$ are normal subgroups with trivial intersection.

  To find the character table of $G/N$ 
  take the submatrix of characters with kernel containing $N$ and delete duplicate columns.
  Since the first row of $Y$ is all ones, this submatrix will be the Kronecker product $[1,1,\dots]\tensor X$.
  Once duplicate columns are removed it will be exactly $X$, and in the same way the character table for $G/H$ is $Y$.

  Now $\idx{G}{H} = \abs{G/H}= \sum_{i} a_{i}^{2}$ and $\idx{G}{N} =\abs{G/N}=\sum_{j} b_{j}^{2}$ so
  \begin{align*}
    \abs{NH} = \abs{N}\abs{H} & = \frac{\abs{G}}{\abs{G/N}} \frac{\abs{G}}{\abs{G/H}}                                                        \\
                              & = \frac{\abs{G}^{2}}{\sum_{j} b_{j}^{2}\sum_{i} a_{i}^{2}} = \frac{\abs{G}^{2}}{\sum_{i,j} (a_{i}b_{j})^{2}} \\
                              & = \frac{\abs{G}^{2}}{\abs{G}} = \abs{G}.
  \end{align*}
  Finally, $G\isomorph N\cross H$ since $N$ and $H$ are normal in $G$, they have trivial intersection, and $G=NH$.
\end{proof}
\end{proposition}

\begin{corollary}
  If $G\cross H$ is transposable then so are $G$ and $H$.
\begin{proof}
  Let $M$ be the character table of $G\cross H$ written as $M=X\tensor Y$ where $X$ and $Y$ are the character tables of $G$ and $H$ respectively.
  Then the (matrix) transpose of $M$ is $X^{T}\tensor Y^{T}$ and by assumption is a character table of some group $\Gamma$ after proper multiplication and division of character degrees.
  It is clear that this multiplication splits across the product so that, by Proposition~\ref{thm:factors}, $\Gamma=G_0\cross H_0$ for some groups $G_0,\, H_0$.
  After appropriate multiplication, $X^{T}$ is the character table of $G_0$ and $Y^{T}$ is that of $H_0$.
  That is, $G$ and $H$ are transposable.
\end{proof}
\end{corollary}

It is well known (see \eg~\cite[Exercise 3.9]{Isaacs:1976}) that the character table determines multiplication of both characters and conjugacy classes.
To be more precise, if $\class_{i}$ is a conjugacy classes and $\classsum_{i}$ its class sum, then class multiplication is
\[ \classsum_{i}\classsum_{j} = \sum_{\nu}a_{ij\nu}\classsum_{\nu} \] %
for non-negative integers $a_{ij\nu}$ which can be calculated by
\[ a_{ij\nu}=\frac{\abs{\class_{i}}\abs{\class_{j}}}{\abs{G}}\sum_{\chi\in\Irr(G)}\frac{\chi(g_{i})\chi(g_{j})\overline{\chi(g_{\nu})}}{\chi(1)}. \] %
Recall that the multiplication constants for the product of two irreducible characters are
\[ b_{ij\nu} = [\chi_{i}\chi_{j},\chi_{\nu}] = \inv{\abs{G}} \sum_{l} \abs{\class_{l}}\chi_{i}(g_{l})\chi_{j}(g_{l})\overline{\chi_{\nu}(g_{l})}. \] %
The following is basically {\cite[Theorem~3.1]{Bannai:1993}} in different terminology.
\begin{proposition}\label{thm:mult-constants}
  Let $G$ and $g^T$ be transposes of each other. Denote by $b_{ij\nu}^{T}$ the multiplication constants for characters of $G^T$ where the indices are inherited from $G$. Then
    \[b_{ij\nu}^{T}= \frac{\sqrt{\classS{\nu}}}{\sqrt{\classS{i}\classS{j}}} \; a_{ij\nu}.\]
\begin{proof}
  Consider conjugacy classes $\class_{i}$ and $\class_{j}$ as characters of $G^{T}$ and compute the coefficients $b_{ij\nu}^{T}$.
  Using the correspondence $\class_{i}\leftrightarrow\phi_{i}$, $\chi\leftrightarrow\class_{\chi}$, and $x_{\chi}\in\class_{\chi}$ as before, we have
  \begin{align*}
    b_{ij\nu}^{T}
    & = \inv{\abs{G}}\sum_{\chi\in\Irr(G)}
    \abs{\class_{\chi}}
    \phi_{i}(x_{\chi})
    \phi_{j}(x_{\chi})
    \overline{\phi_{\nu}(x_{\chi})} \\
    & = \inv{\abs{G}}\sum_{\chi\in\Irr(G)} \chi(1)^{2}
    \frac{\chi(x_{i})\sqrt{\classS{i}}}{\chi(1)}
    \frac{\chi(x_{j})\sqrt{\classS{j}}}{\chi(1)}
    \overline{\parens{\frac{\chi(x_{\nu})\sqrt{\classS{\nu}}}{\chi(1)}}}     \\
    & = \frac{\sqrt{\classS{i}\classS{j}\classS{\nu}}}{\abs{G}}\sum_{\chi\in\Irr(G)}
    \frac{\chi(x_{i})\chi(x_{j})\overline{\chi(x_{\nu})}}{\chi(1)} \\
    & = \frac{\sqrt{\classS{\nu}}}{\sqrt{\classS{i}\classS{j}}}
    \frac{\classS{i}\classS{j}}{\abs{G}}\sum_{\chi\in\Irr(G)}
    \frac{\chi(x_{i})\chi(x_{j})\overline{\chi(x_{\nu})}}{\chi(1)} \\
    & = \frac{\sqrt{\classS{\nu}}}{\sqrt{\classS{i}\classS{j}}} \; a_{ij\nu}.\qedhere
  \end{align*}
\end{proof}
\end{proposition}


\section{Normal Subgroup Correspondences}

In this section we shall investigate normal subgroups of transposable groups, in particular we establish a correspondence between normal subgroups of $G$ and $G^{T}$.
We start with a simple observation.

\begin{proposition}\label{thm:center=abelianization}
  Let $G$ be a transposable group and let $A=G^{T}/(G^{T})'$ be the abelianization of $G^{T}$.
  Then $A\isomorph Z(G)$.
\begin{proof}
  Recall that the linear characters are precisely those whose kernel contains the commutator subgroup.
  Thus we can calculate the character table of $A$ by taking the linear characters of $G^{T}$ and deleting duplicate columns.
  There are exactly enough duplicate columns to make the table square.

  Clearly the linear characters of $G^{T}$ correspond to the central classes of $G$ for these classes have size $1$.
  Note that the number and size of conjugacy classes do not change when we restrict our attention to the character table of $Z=Z(G)$.
  In order to calculate the character table of $Z$, take the character table of $G$ and divide by the character degrees.
  Since all the characters of $G$ are homogeneous on $Z$, each row in the resulting table is an irreducible character of $Z$.

  Next, remove duplicate rows.
  These will be the same as those columns removed when calculating the character table of $A$.
  The rows now correspond to distinct characters and there are enough to make the table square.
  Hence the character table of $Z$ is the transpose of that of $A$.
  Since abelian groups are determined by their character tables we have that $Z\isomorph A$.
\end{proof}
\end{proposition}

In fact we have nearly proven the following more general result.

\begin{proposition}
  Let $G$ be a transposable group and let $N\normal G$ be such that all irreducible characters of $G$ are homogeneous when restricted to $N$.
  Further assume that there is no fusion in $N$, \ie if $n_{1}^{\phan}=n_{2}^{g}$ for $n_{1},n_{2}\in N$, $g\in G$, then $g$ can be chosen to be in $N$.
  Let $\mathcal{N}$ be the set of characters of $G^{T}$ corresponding to the conjugacy classes of $N$.
  Then the character table of $N$ is the transpose of the character table of $H=G^{T}/\ker\mathcal{N}$.
\begin{proof}
  The only thing which does not follow immediately from the previous proof is that the character degrees of $N$ match the square roots of the conjugacy class sizes of $H$.
  This follows from the orthogonality relations and the lack of fusion.
\end{proof}
\end{proposition}


While there is no correspondence which preserves (transposes of) character tables for all $N\normal G$, we can generalize the correspondence to all normal subgroups.
For certain special subgroups, namely the upper and lower central series, we are able to retain some structural information.


\begin{lemma}
  \label{lemma:normal-kernel}
  Given two normal subgroups $N_{1},N_{2}\normal G$ let $N$ be the join $N=N_{1}\join N_{2}=N_{1}N_{2}$ in the normal subgroup lattice.
  It can be determined from the character table by taking the union of the conjugacy classes of $N_{1}$ and $N_{2}$ and finding the set $\mathcal{N}$ of characters whose kernels contain these classes.
  Then the classes in $K=\ker\mathcal{N}$ is the join $N=N_{1}N_{2}$.
\begin{proof}
  Clearly the join $N_{1}N_{2}$ must contain all of the conjugacy classes of both $N_{1}$ and $N_{2}$ (which $K$ does), and it is the minimal such, which $K$ is by construction.
\end{proof}
\end{lemma}

\begin{proposition}
  \label{thm:dual-lattice}
  Let $G$ be a transposable group and $G^{T}$ one of its transpose groups.
  Then for every normal subgroup $N\normal G$, there is a normal subgroup $N^{T}\normal G^{T}$ such that $\abs{G/N}=\abs{N^{T}}$ (equivalently $\abs{N}=\abs{G^{T}/N^{T}}$).
  Furthermore $(N_{1}N_{2})^{T}=N_{1}^{T}\intersect N_{2}^{T}$, so that the lattice of normal subgroups of $G$ is the dual of that of $G^{T}$ including orders of subgroups.
\begin{proof}
  Consider some normal subgroup $N\normal G$ as a collection of conjugacy classes.
  It also can be thought of in terms of the set $\mathcal{N}$ of irreducible characters whose kernels contain $N$.
  In the transposed character table these concepts are switched: the conjugacy classes of $N$ correspond to the characters whose kernel make up $N^{T}$, and the characters defining $N$ become the conjugacy classes of $N^{T}$.

  Consider the normalized character table and note that an entry is in the kernel of a character if its value is $1$.
  In this way kernel entries are easily seen to remain kernel entries after transposition.

  We now determine the orders of $N^{T}$ and $G/N$.
  The character table of $G/N$ is easily determined from the characters in $\mathcal{N}$ by removing duplicate columns.
  In particular this leaves the character degrees intact, so that $\abs{G/N}=\sum_{\chi\in\mathcal{N}}\chi(1)^{2}$.
  The characters in $\mathcal{N}$ become the conjugacy classes of $N^{T}$ whose order can be computed by summing the sizes of its conjugacy classes.
  But this is simply $\abs{N^{T}}=\sum_{\chi\in\mathcal{N}}\chi(1)^{2}$ since the squares of the character degrees in $G$ are the conjugacy class sizes in $G^{T}$.

  To prove that $N_{1}^{T}\intersect N_{2}^{T}=(N_{1}N_{2})^{T}$, observe that conjugacy classes of $N_{1}^{T}\intersect N_{2}^{T}$ are given by the intersection of the set of characters whose kernels are $N_{1}$ and $N_{2}$ respectively.
  The intersection of these sets consists of the irreducible characters whose kernel contains both $N_{1}$ and $N_{2}$, and it follows from Lemma~\ref{lemma:normal-kernel} that $N_{1}^{T}\intersect N_{2}^{T}=(N_{1}N_{2})^{T}$.
\end{proof}
\end{proposition}

We should not be too surprised by this result since characters give information about quotient groups and conjugacy classes about normal subgroups.
Thus interchanging the two concepts might be expected to interchange normal subgroups and quotients as well.



Recall that the lower (descending) and upper (ascending) central series are defined as
\begin{align*}
  \gamma_{1}(G) & = G & \gamma_{i}(G)       & = [\gamma_{i-1}(G),G]; \\
  Z_{0}(G)      & = 1 & Z_{i}(G)/Z_{i-1}(G) & = Z(G/Z_{i-1}(G)).
\end{align*}
One way to determine $Z_{i}(G)$ is to locate the conjugacy classes of $Z_{i-1}(G)$ and then find conjugacy classes $\class$ such that $[\class,G]\subset Z_{i-1}(G)$.
The following proposition gives us a way to do that using only the character table. First define
  \begin{equation}\label{eq:commutator}
    n(g,h)=\sum_{\chi\in\Irr(G)} \frac{\abs{\chi(g)}^{2}\overline{\chi(h)}}{\chi(1)}
  \end{equation}
for $g,h\in G$.

\begin{proposition}%
  \label{thm:commutator}
Let $G$ be a group. For $g\in G$ we have
\begin{equation}\label{eq:Zi}
  g\in Z_i(G)\Leftrightarrow n(g,h)=0\,\forall h\not\in Z_{i-1}(G).
\end{equation}
\begin{proof}
  By \cite[Problem 3.10]{Isaacs:1976} $n(g,h)\ne0$ if and only if $h$ is conjugate to $[g,g']$ for some $g'\in G$.

  If $g\in Z_i(G)$ then this cannot be satisfied by $h\not\in Z_{i-1}(G)$.

  Conversely, if $n(g,h)\ne 0$ implies $h\in Z_{i-1}$ then $[g,G]\subseteq Z_{i-1}$ so $g\in Z_i$.
\end{proof}
\end{proposition}

We now translate this description into a condition within $G^{T}$.
Let $\class_{g}$ denote the $G$-conjugacy class of $g$, and $\phi_{g}$ the character of $G^{T}$ which corresponds to this class.
If $\chi$ is an irreducible character of $G$, then $\class_{\chi}$ will denote the corresponding conjugacy class of $G^{T}$, and $x_{\chi}$ a representative of this class.
With this notation (as the normalized character tables are transpose) we have the identity
\begin{equation*}
  \frac{\chi(g)}{\chi(1)} = \frac{\phi_{g}(x_{\chi})}{\phi_{g}(1)}.
\end{equation*}
Substituting into $n(g,h)\ne0$ in~\eqref{eq:commutator} gives
\begin{align*}
  \sum_{\class_{\chi}} \frac{ \abs{\phi_{g}(x_{\chi})}^{2} \frac{\chi(1)^{2}}{\phi_{g}(1)^{2}}  \overline{\phi_{h}(x_{\chi})} \frac{\chi(1)}{\phi_{h}(1)}}
  {\chi(1)}                                         & \not= 0 \\
  \sum_{\class_{\chi}} \frac{ \abs{\phi_{g}(x_{\chi})}^{2}\overline{\phi_{h}(x_{\chi})} {\chi(1)^{2}} }
  {\phi_{g}(1)^{2}\phi_{h}(1)}                      & \not= 0 \\
  \intertext{and since $\phi_{g}(1)^{2}\phi_{h}(1)$ is constant}
  \sum_{\class_{\chi}} \abs{\phi_{g}(x_{\chi})}^{2}\overline{\phi_{h}(x_{\chi})} {\abs{\class_{\chi}}}
                                                    & \not= 0 \\
  \abs{G^{T}}[\phi_{g}\overline{\phi_{g}},\phi_{h}] & \not= 0 \\
  [\phi_{g}\overline{\phi_{g}},\phi_{h}]            & \not= 0
\end{align*}
where $[\ ,\ ]$ is the usual inner product on characters of $G^{T}$.

Now denote $n^T(\phi,\psi)=[\phi\overline{\phi},\psi]$.

\begin{proposition}\label{thm:char-central-series}
  Let $G$ be a group. For $\phi\in\Irr(G)$ we have
\begin{equation}\label{eq:gammai}
  \ker\phi\supseteq\gamma_i(G)\Leftrightarrow n^T(\phi,\psi)=0\, \forall \psi,\ker\psi\nsupseteq \gamma_{i-1}(G).
\end{equation}
\begin{proof}
  Let $U$ be a module affording $\phi$, so that $\phi\overline{\phi}$ is the character of $\Hom(U,U)$.
  We need to prove that the kernel of $U$ contains $\gamma_{i}(G)$ if and only if every simple submodule of $\Hom(U,U)$ has kernel containing $\gamma_{i-1}(G)$, that is if and only if $\Hom(U,U)$ has kernel containing $\gamma_{i-1}(G)$.

  If $\gamma_{i}(G)\subset\ker U$ then every $c\in\gamma_{i-1}(G)$ is central with respect to the action on $U$.
  Let $\mu\in\Hom(U,U)$, and $c\in\gamma_{i-1}(G)$, then we wish to prove that $c$ fixes $\mu$, \ie $\mu c=\mu$.
  Consider the action on a vector $u\in U$, by definition
  \begin{align*}
    (\mu c)(u) & = \mu(uc^{-1})c.            \\
    \intertext{Now, the action of $c$ is central on $U$, so it is equivalent to multiplication by a scalar $\lambda$}
               & = \mu(u\lambda^{-1})\lambda \\
               & = \mu(u).
  \end{align*}
  This is true for every $\mu$ and every $u$, so that $\ker \Hom(U,U) \supseteq \gamma_{i-1}(G)$.

  To prove the converse, suppose $\ker \Hom(U,U) \supseteq \gamma_{i-1}(G)$.
  This means that for $c\in\gamma_{i-1}(G)$ we have $\mu c=\mu$, that is
  \begin{align*}
    \mu(u)       & = (\mu c) (u) = \mu (uc^{-1})c.
    \shortintertext{Then }
    \mu(u)c^{-1} & =\mu(uc^{-1})
  \end{align*}
  for all $\mu$ and all $u$, so the action of $c$ commutes with every $\mu$, hence the action of $c$ is central in $U$.
  This is true for every $c$ in $\gamma_{i-1}(G)$, so the kernel of $U$ contains $[\gamma_{i-1}(G),G]=\gamma_{i}(G)$ as claimed.
\end{proof}
\end{proposition}

\begin{theorem}
  Given a transposable group $G$, the (abelian) factors $Z_{i}(G)/Z_{i-1}(G)$ are isomorphic to $\gamma_{i}(G^{T})/\gamma_{i+1}(G^{T})$ for all $i$.
\begin{proof}
  From Propositions~\ref{thm:commutator}~and~\ref{thm:char-central-series} there is a correspondence between the conjugacy classes of the upper central series of $G$ and the irreducible characters defining the lower central series of $G^{T}$.
  Due to the special nature of the subgroups involved we are able to determine the character tables (and hence isomorphism type) of the central series factors, from the character table of $G$.

  To find the character table of $Z_{i}(G)/Z_{i-1}(G)$, restrict attention to the characters of $G$ which contain $Z_{i-1}(G)$ in their kernel, and the conjugacy classes of $Z_{i}(G)$.
  This corresponds to finding the character table of $G/Z_{i-1}(G)$, and then restricting to $Z_{i}(G)/Z_{i-1}(G)$.
  As in Proposition~\ref{thm:center=abelianization}, we simply divide by character degrees and then remove duplicate columns and duplicate rows to find the character table of $Z_{i}(G)/Z_{i-1}(G)$.

  Finding the character table of $\gamma_{i}(G^{T})/\gamma_{i+1}(G^{T})$ is done in a completely dual way.
  This time we restrict attention to the conjugacy classes of $\gamma_{i}(G^{T})$, and characters whose kernel contains $\gamma_{i+1}(G^{T})$.
  These are the same rows as columns before, and vice versa.
  As before, we divide by character degrees and then remove duplicate rows and columns, leading to the transpose of the character table of $Z_{i}(G)/Z_{i-1}(G)$.
  Since they are abelian groups, $Z_{i}(G)/Z_{i-1}(G)$ and $\gamma_{i}(G^{T})/\gamma_{i+1}(G^{T})$ are isomorphic.
\end{proof}
\end{theorem}

\section{Nilpotency of Transposable Groups}

This section makes use of the theory of $p$-blocks.
We refer the reader to~\cite{Isaacs:1976} for basic definitions and results.

Let $\varepsilon$ be a primitive $|G|$-th root of unity and $R$ the algebraic integers of $\mathbb{Q}[\varepsilon]$.
Throughout this section we will write $a\equiv b$ to mean that $a$ is congruent to $b$ modulo some fixed maximal ideal of $R$ containing the fixed prime $p$.

\begin{proposition}\label{thm:full-defect}
  Every $p$-block $B$ of a transposable group $G$ is of full defect.
\begin{proof}
  By the definition of defect of a block $B$, there is $\chi\in B$, and a conjugacy class $\class$ (a defect class) such that
  \begin{equation*}
    \chi(x^{-1})\frac{\abs{\class}\chi(x)}{\chi(1)}\not\equiv0
  \end{equation*}
  for $x\in\class$.
  By the definition of transposable groups (and using the correspondence $\class\leftrightarrow\phi_{\class}$, and $\chi\leftrightarrow\class_{\chi}$ with $x_{\chi}\in\class_{\chi}$), we have that
  \begin{equation*}
    \frac{\abs{\class}\chi(x)}{\chi(1)} = \phi_{\class}(1)\phi_{\class}(x_{\chi})\not\equiv0
  \end{equation*}
  and in particular $\phi_{\class}(1)\not\equiv0$, so that $\abs{\class}=\phi_{\class}(1)^{2}\not\equiv0$.
  Therefore, $B$ has full defect.
\end{proof}
\end{proposition}

\begin{proposition}\label{thm:KK}
  Let $G$ be transposable.
  If $\chi\in B_{0}(G)$, the principal block of $G$ then for all $\phi\in\Irr(G^{T})$, $\frac{\abs{\class_{\chi}}\phi(x_{\chi})}{\phi(1)}\equiv\abs{\class_{\chi}}$.
\begin{proof}
  Let $\phi\Irr(B)$ for a block $B$ of $G^T$. Since
  \begin{align*}
    \frac{\abs{\class_{\chi}}\phi_{1}(x_{\chi})}{\phi_{1}(1)}\equiv\frac{\abs{\class_{\chi}}\phi_{2}(x_{\chi})}{\phi_{2}(1)}
  \end{align*}
  for all $\phi_{1}$ and $\phi_{2}$ in $B$, we may assume that $\phi(1)\not\equiv0$ by Proposition~\ref{thm:full-defect}. Let $\phi$ correspond to $\class$, so $\phi_{\class}=\phi$.

  Now $\chi\in B_{0}(G)$, so for $\class$
  \begin{align*}
    \frac{\abs{\class}\chi(x)}{\chi(1)}     & \equiv \abs{\class} \\
    \phi_{\class}(1)\phi_{\class}(x_{\chi}) & \equiv \phi_{\class}(1)^{2},
    \shortintertext{so by $\phi_{\class}(1)\not\equiv1$ we get}
    \frac{\phi_{\class}(x_{\chi})}{\phi_{\class}(1)}& \equiv 1.
  \end{align*}
  From this follows the claimed $\frac{\abs{\class_{\chi}}\phi(x_{\chi})}{\phi(1)}\equiv\abs{\class_{\chi}}$.
\end{proof}
\end{proposition}

\begin{corollary}\label{coro:B0-divis-by-p}
  If $G$ is transposable and $\chi\in B_{0}(G)$ with $\chi(1)\not\equiv0$, then $\chi(1)=1$.
\begin{proof}
  Since $\abs{\class_{\chi}}=\chi(1)^{2}\not\equiv0$, it follows from Proposition~\ref{thm:KK} that $\chi(1)\chi(x)=\frac{\abs{\class_{\chi}}\phi(x_{\chi})}{\phi(1)}\not\equiv0$ for all $x\in G$.
  In particular $\chi(1)\chi(x)\not=0$ for all $x\in G$, and so $\chi$ must be linear by a famous theorem of Burnside~\cite[3.15]{Isaacs:1976}.
\end{proof}
\end{corollary}

\begin{theorem}[{\cite[Corollary~3]{Isaacs:1976a}}]\label{thm:B0-p-nilpotent}
  Suppose that every nonlinear irreducible character of $G$ in the principal $p$-block has degree divisible by $p$.
  Then, $G$ has a normal $p$-complement, \ie is $p$-nilpotent.
\end{theorem}

\begin{corollary}\label{coro:transposable=nilpotent}
  If $G$ is transposable then $G$ is nilpotent.
\end{corollary}




\section{Examples of Transposable Groups}
\label{sec:examples-trans-gps}


In this section we consider non-abelian transposable groups.
Since such groups are nilpotent, and transposability is determined by direct factors, we only consider directly indecomposable $p$-groups.

We first give a method of constructing transposable groups by combining transposable groups.
This method yields abelian groups when combining abelian groups, but allows us to turn a nonabelian transposable group into a family of such groups.
This allows us to concentrate on certain basic families of groups which we call stem groups, mimicking the terminology of isoclinism.
When adding cyclic groups, our method is reminiscent of growing roots or branches onto a group so the terminology is fitting.

After we discuss this construction we will give examples of all known stem groups.
The stem groups are self-dual, but our construction allows us to create groups which are not.
It is an interesting question whether there are any stem groups that are transposable but not self-dual.


\begin{proposition}
  Let $G$ be a transposable group, $x\in Z(G)$ and $X=\angles{x}$.
  Furthermore suppose that $X\intersect G'=1$.
  Then $G/X$ is transposable.
\begin{proof}
  By induction on $|X|$ we can assume $o(x)=p$ a prime.

  Take $y\in G\setminus X$ arbitrary. Then $\{(yX)^{tX}\mid tX\in G/X\}=\{y^tX\mid tX\in G/X\}$ has the same size as the conjugacy class of $y\in G$. It applies equally to $y,yx,yx^2,\ldots,yx^{p-1}$. The character table of $G/X$ is determined by taking the characters which contain $x$ in their kernel and then deleting duplicate columns. That is, $p-1$ columns are deleted for each one remaining. In particular the character values do not change.




  Let now $\lambda\in\Irr(G^T)$ correspond to the class $\{x\}$ of $G$. Hence, $H=\ker\lambda$ has index $p$ in $G^{T}$. Since $x\not\in G'$, the kernel of $\lambda$ does not contain the center of $G^{T}$.
  Thus, $Z(G^{T})$ contains a full set of representatives of right-cosets of $H$.
  Therefore
  $\phi_{H}$ is homogeneous for all $\phi\in\Irr(G)$.
  But $H$ being of index $p$ it means $\phi_{H}$ is irreducible. (See for example \cite[6.19]{Isaacs:1976}.) Each such $\phi_{H}$ is a restriction of $p$ distinct irreducibles of $G^T$
  And now the character table of $H$ is determined by restricting the rows to the classes of $H$
  and then deleting duplicate rows. That is, $p-1$ rows are deleted for each one remaining. In particular the character values do not change and we obtain the transpose of that of $G/X$.
\end{proof}
\end{proposition}

\begin{corollary}
  If $G$ is transposable and $N\normal G$ such that $G=NZ(G)$, then $N$ is transposable.
\end{corollary}

By starting with two transposable groups $G_{1},G_{2}$ we can create new (directly indecomposable) transposable groups as factors or subgroups of $G_{1}\cross G_{2}$.
If $G_{2}$ is abelian, then we can think of it as growing roots or branches.
After growing such roots or branches, the commutator structure does not change, rather only the power map.
All of the known examples can be constructed in this way from self-dual groups.

In particular, starting with a transposable group $G$ and $g,h\in G$ with $g^{p}=h^{p}=1$, $g\not\in G'$, $h\in Z(G)$ we can create a new group with a new power map of $g^{p}=h$ and $h^{p}=1$.


\subsection{Stem Groups}

Apart from abelian groups there are two families of finite transposable groups found in the literature, to which we add one more.

The first family, discussed in~\cite{Hanaki:1997a}, exists for $p\geq3$ and has order $p^{5}$.
It has a presentation of the form
\begin{align*}
  G = \presentation{a_{1},a_{2},b,c_{1},c_{2}}{[a_{1},a_{2}]=b,[a_{i},b]=c_{i},a_{i}^{p}=\zeta_{i},b_{\phan}^{p}=c_{i}^{p}=1}
\end{align*}
where $\zeta_{i}$ is central.
The values of $\zeta_{i}$ leading to distinct groups is found in~\cite{James:1980}, but they are irrelevant for our purposes.
In fact, by growing roots and branches appropriately, we can get them all from the case when $\zeta_{i}=1$.

It is clear from the presentation that $Z(G)=\angles{c_{1},c_{2}}\cong \mathbb{Z}_{p}^{2}$ and $G'=\angles{b,c_{1},c_{2}}$.
Moreover, $G/Z(G)$ is an extraspecial group of type $+$.
The nilpotency class is $3$, and derived length is $2$.
All normal subgroups of $G$ are either contained in the center, or contain the commutator subgroup, which leads to the normal subgroup lattice in Figure~\ref{fig:nsl-class3}.
The number of normal subgroups in the center is $p+1$.
\begin{figure}[h]
  \begin{minipage}[b]{.5\linewidth}
    \centering
    \begin{tikzpicture}[scale=1.5]
      \tikzstyle{every node}=[rectangle, fill=white, inner sep=2pt, minimum width=3pt]
      \foreach \y in {0,3} {
        \foreach \x in {-1,-0.6,-0.2,0.2,1} {
          \draw (0,\y) -- (\x,\y+1) -- (0,\y+2);
        }
        \draw (0.6,\y+1) node{$\cdots$};
      }
      \draw (0,0) node{$\phantom{1}1\phantom{1}$}
      (0,2) node{$Z(G)$}
      -- (0,3) node{$G'$}
      (0,5) node{$G$};
    \end{tikzpicture}
    \vspace*{1.5cm}
    \caption{Hanaki's Groups}\label{fig:nsl-class3}
  \end{minipage}%
  \begin{minipage}[b]{.5\linewidth}
    \centering
    \begin{tikzpicture}[scale=1.5]
      \tikzstyle{every node}=[rectangle, fill=white, inner sep=2pt, minimum width=3pt]
      \foreach \y in {0,5} {
        \foreach \x in {-1,-0.6,-0.2,0.2,1} {
          \draw (0,\y) -- (\x,\y+1) -- (0,\y+2);
        }
        \draw (0.6,\y+1) node{$\cdots$};
      }
      \draw (0,0) node{$\phantom{1}1\phantom{1}$}
      (0,2) node{$Z(G)$}
      -- (0,3) node{$Z_{2}(G)$}
      -- (0,4) node{$Z_{3}(G)$}
      -- (0,5) node{$G'$}
      (0,7) node{$G$};
    \end{tikzpicture}
    \caption{Class 5 Groups}\label{fig:nsl-class5}
  \end{minipage}
\end{figure}

The second family appears to be unknown in the literature as being self-dual.
It exists for $p\geq5$ and has order $p^{7}$.
Its structure is 
similar to the previous family.
For the cases with $p\leq11$, these groups are included in the small groups database available in GAP~\cite{GAP4}.
It has a presentation of
\begin{align*}
  G = \big\langle a_{1},a_{2},b_{1},b_{2},b_{3},c_{1},c_{2} \big\vert
  & [a_{1},a_{2}]=b_{1}, [a_{1},b_{i=1,2}]=b_{i+1}, [a_{1},b_{3}]=c_{1},\\
  & [a_{2},b_{2}] = [a_{2},b_{3}] = [b_{2},b_{1}] = c_{2}, \\
  &  a_{i}^{p}=\zeta_{i},b_{i}^{p}=c_{i}^{p}=1\big\rangle
\end{align*}
where $\zeta_{i}$ is central.

It is clear that $Z(G)=\angles{c_{i}}\cong \mathbb{Z}_{p}^{2}$ and $G'=\angles{b_{i},c_{i}}$.
The nilpotency class is $5$, and the derived length is $2$.
Similar to before, $G$ has coclass $2$ and $G/Z_{3}(G)$ is extraspecial group of type $+$.
Its normal subgroup lattice is also similar and given in Figure~\ref{fig:nsl-class5}.

The proof that this family of groups is transposable will appear in the Ph.D.~dissertation of the first author.

The final family of self-dual groups is a generalization of the Suzuki $2$-groups $A(n,\theta)$, and has been studied in several papers~\cite{Riedl:1999,Hanaki:1996a,Sagirov:2003}.
The non-trivial claims below are proven in one or more of those papers.

Let $q=p^{a}$ be a prime power, $s$ and $l$ be positive integers and let $\theta$ be a generator of the Galois group of $\mathbb{F}_{q^{s}}$ over $\mathbb{F}_{q}$.
Furthermore, assume that $s$ is odd and $(s,l!)=(s,q)=(s,q-1)=1$.
Then we can define a group, $G=G(q,s,l)$, whose elements are $l$-tuples $a=(a_{1},\dots,a_{l})$ of elements of $\mathbb{F}_{q^{s}}$ and multiplication $c=ab$ given by
\begin{equation*}
  c_{i} = a_{i}+ \sum_{j=1}^{i-1}a_{i-j}^{\theta^{j}}b_{j}^{\phan} +b_{i}.
\end{equation*}
This can be thought of as multiplication of skew polynomials of the form $1+\sum_{i=1}^{l}a_{i}x^{i}$ modulo the ideal $\parens{x^{l+1}}$, see~\cite{Riedl:1999}.
It is also equivalent to
\begin{align*}
  a=(a_{1},a_{2},\dots,a_{l})
  =
  \begin{pmatrix}
    1                                                                  \\
    a_{1}  & 1                                                         \\
    a_{2}  & a_{1}^{\theta}   & 1                                      \\
    a_{3}  & a_{2}^{\theta}   & a_{1}^{\theta^{2}}   & 1               \\
    \vdots & \vdots           &                      & \ddots & \ddots \\
    a_{l}  & a_{l-1}^{\theta} & a_{l-2}^{\theta^{2}} & \cdots & a_{1}^{\theta^{l-1}} & 1
  \end{pmatrix}
\end{align*}
with the regular matrix multiplication.

The lower and upper central series of such a group $G$ coincide and each central series factor is isomorphic to the additive group of $\mathbb{F}_{q^{s}}$ \ie elementary abelian of order $q^{s}$.
In fact, the central series is described in terms of ``layers'' in the matrix or fixed powers of $x$.
The nilpotency class of $G$ is $l$, and its derived length is $\lceil\log_{2}(l+1)\rceil$.
Clearly $G(q,s,l)/Z(G(q,s,l))\cong G(q,s,l-1)$.

The group $G$ has $q^{l-i}(q^{s}-1)$ (or $q^{s}$ if $l=i$) irreducible characters of degree $q^{(l-i)(s-1)/2}$ and the same number of conjugacy classes of size $q^{(l-i)(s-1)}$.

Not all such groups are self-dual, however.
If $l>p$ then $G(p^{a},s,l)$ is not self-dual, but when $l=s-1<p$ it is~\cite{Hanaki:1997}.
It is unknown whether $G(p^{a},s,l)$ is self-dual for all $l\leq p$, but at least $G(p,3,2)$ is self-dual for all $p$, so it is reasonable to conjecture that they are.



\bibliographystyle{alpha}
\bibliography{transposable4}

\end{document}